\documentclass{article}
\usepackage{spconf,amsmath,graphicx,amssymb,mathtools}
\usepackage{xcolor}
\usepackage{subcaption}
\usepackage[free-standing-units, space-before-unit, use-xspace]{siunitx}

\newcommand{\auto}{\mathbf N}
\newcommand{\encoder}{\mathbf E}
\newcommand{\decoder}{\mathbf D}

\newcommand{\signal}{u}
\newcommand{\signalP}{\signal_\Plus}
\newcommand{\data}{y}
\newcommand{\Ko}{\mathbf K}
\newcommand{\Ro}{\mathbf{K}^\sharp}
\newcommand{\Uo}{\mathbf U}
\newcommand{\X}{\mathbb X}
\newcommand{\Y}{\mathbb Y}

\newcommand{\M}{\mathcal M}

\newcommand*{\N}{\mathbb{N}}
\newcommand{\B}{\boldsymbol{\Delta}}

\def\Plus{\texttt{+}}
\newcommand{\reg}{\mathcal{R}}    
\newcommand{\tik}{\mathcal{T}}     

\DeclareMathOperator*{\argmin}{arg\,min}

\DeclarePairedDelimiter{\abs}{\lvert}{\rvert}
\DeclarePairedDelimiter{\norm}{\lVert}{\rVert}

\newcommand{\kl}[1]{\left(#1\right)}

\newtheorem{theorem}{Theorem}

\numberwithin{equation}{section}
\numberwithin{figure}{section}
\numberwithin{equation}{section}
\numberwithin{table}{section}
\numberwithin{theorem}{section}

\title{SPARSE aNETT FOR SOLVING INVERSE PROBLEMS WITH DEEP LEARNING}
\name{Daniel Obmann$^1$, Linh Nguyen$^2$, Johannes Schwab$^1$, Markus Haltmeier$^{1,}$\sthanks{Corresponding author: { markus.haltmeier@uibk.ac.at}. DO and MH  acknowledge support of the Austrian Science Fund (FWF), P 30747-N32.}$^,$\sthanks{Due to page restrictions the official proceeding (ISBI 2020) only contains 4 pages.}}
\address{$^1$University of Innsbruck, Department of Mathematics, Innsbruck, Austria\\
	$^2$University of Idaho, Department of Mathematics, Moscow, US}

\begin{document}
%
\maketitle
\begin{abstract}
We propose a  sparse reconstruction framework  (aNETT) for solving inverse problems. Opposed to existing sparse reconstruction techniques that are  based  on  linear sparsifying transforms,  we train an  autoencoder network $\decoder \circ \encoder$ with $\encoder$ acting as a nonlinear sparsifying transform and minimize a Tikhonov functional with learned regularizer formed by the $\ell^q$-norm of the encoder coefficients and  a penalty for the distance to the data manifold. We propose a strategy for training an autoencoder based on a sample set of the  underlying image class such that the autoencoder is independent of the forward operator and is subsequently adapted to the specific forward model. Numerical results are presented for sparse view CT, which clearly demonstrate the  feasibility, robustness and the  improved  generalization capability and stability of aNETT over post-processing networks. 
\end{abstract}
\section{Introduction}

Various  applications in medical imaging, remote sensing and
elsewhere require solving  inverse problems of the form
\begin{equation}\label{eq:ip}
\data   = \Ko \signal + z \,,
\end{equation}
where $\Ko \colon \X \to \Y$ is a linear operator between
Hilbert spaces, and $z$ is the  data distortion.
Inverse problems are well analyzed and several established
approaches for its solution exist \cite{EngHanNeu96,scherzer2009variational}.
Recently, neural networks (NN) and deep learning  appeared as a new paradigms for solving inverse problems and demonstrate  impressive  performance \cite{lee2017deep,jin2017deep,sun2016deep,wang2016perspective}.

In order to enforce data consistency, in \cite{li2020nett} a  deep  learning
approach named NETT (NETwork Tikhonov Regularization) has  been proposed and analyzed
based on minimizing $\norm{\Ko(\signal) - y}^2  +  \alpha \reg_0(\signal )
$, where $ \reg_0 $ is a trained network serving as regularizer.
One of the main assumptions for the analysis of \cite{li2020nett}  is the coercivity of the regularizer which requires special care in network design and training. 
In order to overcome this limitation, we introduce the sparse augmented NETT (aNETT), which considers minimizers of
\begin{multline}\label{eq:anett}
\tik_{\alpha,  y}( \signal )
\coloneqq
\norm{ \Ko \signal - y}^2
\\ +
\alpha  \kl{  \sum_{\lambda \in \Lambda } w_\lambda \abs{(\encoder (\signal))_\lambda}^q
 + \frac{c}{2} \norm{\signal - \auto ( \signal )}^2 } \,.
\end{multline}
Here $\auto = \decoder \circ \encoder$ is a sparse  autoencoder network,
$\encoder  \colon  \X \to \ell^2(\Lambda)$ and $\decoder  \colon  \ell^2(\Lambda) \to \X $
are the encoder and decoder network, 
$\Lambda$ is a countable index set, and 
 $\ell^2(\Lambda)$ is the latent Hilbert space of sparse codes.  
The  weighted $\ell^q$-norm  $\norm{\encoder(\signal)}_{q,w} \triangleq \sum_{\lambda \in \Lambda } w_\lambda \abs{(\encoder (\signal))_\lambda}^q$ implements learned sparsity, and the augmented term $\norm{\signal - \auto ( \signal) }^2$ is to force $\signal$ to be close to the data manifold $\M$. Both terms together allow to show coercivity of the regularizer. 
Based on this we derive  stability, convergence and convergence rates  for aNETT. Note that sparse regularization is well investigated for linear representations   \cite{GraHalSch08,daubechies2004sparsity} but so far has not been investigated for nonlinear deep autoencoders.

\section{Sparse augmented NETT}

\subsection{Theoretical results} \label{sec:analyticresults}

Throughout this section we assume the following.
\begin{itemize}
  \setlength\itemsep{-0.5em}
  \item $\Ko \colon \X \to \Y$ is linear and bounded.
  \item $\encoder \colon \X \to \ell^2(\Lambda)$ is weakly sequentially continuous.
  \item $\decoder \colon \ell^2(\Lambda) \to \X$ is weakly sequentially continuous.
  \item $q \geq 1$, $c >0$, 
  \item $w_{\rm min} \triangleq \inf \{w_\lambda \mid \lambda \in \Lambda \} >0$.
\end{itemize}
Furthermore, we define $\auto \triangleq \decoder \circ \encoder$ and choose the regularizer
\begin{equation*}
\reg_c(\signal)  \triangleq \sum_{\lambda \in \Lambda } w_\lambda \abs{(\encoder (\signal))_\lambda}^q+\frac{c}{2} \norm{\signal - \auto( \signal )}^2.
\end{equation*}
Under these assumptions, \eqref{eq:anett}
has a minimizer for all $y \in \Y$ and all $ \alpha>0$.
Moreover, we have the following results.

\begin{theorem}[Convergence]\label{thm:conv}
Let   $\data \in \Ko(\X)$, $\data_k\in \Y$ for $k\in\N$ satisfy $\norm{\data_k -\data} \leq  \delta_k$, and $\delta_k, \alpha_k, \delta_k^2 / \alpha_k \to 0$
as $k \to \infty$.
 Then with $\signal_k \in \argmin_\signal \tik_{\alpha_k, \data_k} (\signal)$  the following hold:
\begin{itemize}
\setlength\itemsep{-0.5em}
\item $(\signal_k)_{k \in \N}$ has at least one weak  accumulation point.

\item Every weak accumulation point of $(\signal_k)_{k\in \N}$ is an
$\reg_c  $-minimizing solution of $\Ko\signal = \data$.

\item If $\Ko \signal = \data$ has a unique $\reg_c $-minimizing solution $\signalP$, then $(\signal_k)_{k \in \N}$  weakly converges to $\signalP$.
\end{itemize}
\end{theorem}

\begin{theorem}[Convergence rate]\label{thm:rates}
Let $\reg_c$ be G\^ateaux differentiable, 
$\Ko$ have finite-dimensional range
and consider minimizers $\signal_{\alpha,\delta} \in \argmin_\signal \tik_{\alpha, \data_\delta} (\signal)$ with $\norm{\data_\delta-\Ko \signal_\Plus} \leq \delta$.
Then $\alpha \sim \delta$ implies the convergence rate $\B_c(\signal_{\alpha,\delta}, \signal_\Plus)   = \mathcal{O}(\delta)$ as $\delta \to 0$ in terms of the so-called absolute Bregman distance
$\B_c(\signal, \signal_\Plus) \triangleq \abs{  \reg_c(\signal) - \reg_c(\signal_\Plus) - \langle \nabla \reg_c(\signal_\Plus), \signal - \signal_\Plus \rangle } $.
\end{theorem}

Proofs of Theorems \ref{thm:conv} and \ref{thm:rates} are given in~\cite{haltmeier2019sparse}.

\subsection{Trained autoencoder}
\label{ssec:train}

First, an  autoencoder $\auto^{a}$ is trained such that  $\signal$ is close to  $\auto^{a} ( \signal) $ and  that $\norm{\encoder(\signal)}_{1,w}$ is small for any $\signal$ in a class $\M$ of images of interest. For that purpose, we add the regularizer $\norm{\,\cdot\,}_{1,w}$ to the loss function for training $\auto^{a}$ as denoising network. To be more specific, let  $(\auto_\theta^{a})_{\theta \in \Theta}$ be a  family of autoencoder networks $\auto_\theta^{a} = \decoder_\theta^{a} \circ \encoder_\theta$, where $\encoder_\theta \colon \X \to \ell^2(\Lambda)$ are admissible (in the sense of above assumptions) encoder  networks and   
$\decoder_\theta^{a} \colon \ell^2(\Lambda) \to \X$  admissible decoder networks.  Moreover, suppose that 
$\signal_1, \dots, \signal_m \in \M$ is a training dataset. To select the particular autoencoder based on the training data, we  consider the following  training strategy for the sparse denoising autoencoder
\begin{multline}\label{eq:AEtrain}
\theta^*  \in \argmin_{\theta} \frac{1}{m}\sum_{i=1}^m \norm{\auto_\theta^{a} (\signal_i + \varepsilon_i) - \signal_i}_2^2 \\
+ \eta \norm{\encoder_\theta(\signal_i + \varepsilon_i) }_{1, w} + \beta  \norm{\theta}_2^2  \,,
\end{multline}
and set  $[\auto^{a}, \encoder] \triangleq [\auto_{\theta^*}^{a}, \encoder_{\theta^*}]$. Here  $v_i \in \X$ are data perturbations and  $\beta > 0$  a regularization parameter. 

By training with perturbed data points $\signal_i + \varepsilon_i$, we increase  robustness of the trained autoencoder.  Note that  the perturbations $\varepsilon_i \in \X$ are chosen independently of the operator $\Ko$ such that the autoencoder can be used  for each forward operator in a universal manner. Clearly then, the autoencoder depends on the specific  manifold $\M$ of images of interest.  As we shall see however,  opposed to typical  deep learning based reconstruction methods which do not account for data consistency outside the training data set, the sparse aNETT is  robust against changes of the specific image manifold.       
Note that Thms. \ref{thm:conv} and \ref{thm:rates}
hold true for   $\auto^{a}$ in place of $\auto$.

\subsection{Adaptation to specific forward models}

The sparse aNETT  \eqref{eq:anett} consists of a 
data consistency term,
a  sparsity term, and   an augmented  term  enforcing  
$\auto(\signal) \simeq  \signal$.
Ideally,  the set of all  approximately data consistent elements
that are also approximate fixed points of $\auto$, is close to the image manifold $\M$.  However, without adjusting the autoencoder to specific forward
 models, this is a challenging and maybe impossible task. 
 Indeed, for the application we consider in this paper, namely sparse view CT, we observed that the autoencoder trained independent of the forward operator, was not able to sufficiently well distinguish  between data-consistent elements inside and outside desired image class.

One way to increase the value of $\norm{\signal - \auto( \signal )}$ for undesired but data consistent   elements is to adopt the training strategy developed in \cite{li2020nett}    and to take  the  data perturbations in \eqref{eq:AEtrain} as  $\varepsilon_i = \Ro \Ko \signal_i - \signal_i$ where  
$\Ro$ is a reconstruction operator approximating the  Moore-Penrose inverse of $\Ko$, $\signal_i$ are the artifact free images and $ \Ro \Ko \signal_i $  images with artefacts.   
In this case, the training dataset depends on the forward operator, and the autoencoder has to be retrained for every specific forward operator.  Therefore, in this paper we follow a different approach. Instead of adjusting the autoencoder training,  we compose the operator independent  autoencoder $\auto^{a}$ with another network $\Uo$, that  is trained to distinguish between the desired  images and images  with operator dependent artefacts. For that purpose we choose a network architecture  $( \Uo_\kappa)_{\kappa \in K}$ and  select  $\Uo = \Uo_{\kappa^*}$,  where $\kappa^*$ is a minimizer of
\begin{align}\label{eq:U} 
\frac{1}{2m}\sum_{i=1}^{2m} \norm{\Uo_\kappa (\auto^{a}(v_i)) - \signal_i}_2^2 + \gamma \norm{\kappa}_2^2 \,,
\end{align}
where $v_i = \Ro \Ko \signal_i$ for $i = 1, \dots, m$ and $v_i = \signal_i$ for $i = m+1, \dots, 2m$ and $\gamma > 0$ is a regularization parameter. We see  that Thms. \ref{thm:conv} and \ref{thm:rates}
still hold true  for the final autoencoder  $\auto \triangleq \Uo \circ \auto^{a}$ if $\Uo$ is weakly sequentially continuous.


\section{Application to sparse view CT} 
\label{sec:ct}

For the numerical simulations we consider the problem of recovering an image from sparse view parallel-beam CT data with $60$ angles. For this problem, the forward operator $\Ko$ is given by the angularly  subsampled Radon transform
\begin{align*}
(\Ko u)(s, \varphi ) \triangleq  \int_{L(s, \varphi)} \signal (x) \mathrm{d}\sigma(x) \,,
\end{align*}
for $60$ equidistant angles $\varphi$  in  $[0, \pi]$.
Here $L(s, \varphi)$ is the line in the plane with normal vector $(\cos(\varphi), \sin(\varphi) )$ and signed distance $s \in [-1.5,1.5]$ from the origin.   
Discretization of the Radon transform  is done using the ODL library \cite{adler2017odl}. 
The data chosen for the numerical simulations are taken from the  Low Dose CT Grand Challenge \cite{mccollough2016tu}. We consider the images at \SI{1}{\milli\metre} slice thickness given in the dataset and take the first seven patients for training (4267 images), the next two patients for validation (1143 images) and the last patient for testing (526 images). Each of these images is rescaled to have pixel values in the interval $[0,1]$.

\subsection{Network training}

We first train  $ \auto^{a}$, $\encoder $ by minimizing  \eqref{eq:AEtrain} and  subsequently  train $\Uo$ by minimizing \eqref{eq:U}.   The sparse autoencoder is chosen as $\auto = \Uo \circ  \auto^{a}  $.
The network architecture chosen for the problem adapted network $\Uo$ is the tight frame U-Net \cite{han2018framing} and the auto-encoder architecture is  chosen as in \cite{obmann2020deep}. 
The perturbations in \eqref{eq:AEtrain} are taken as independent realizations of Gaussian white noise with  standard deviation $p \cdot \bar{\signal}_i$ where $p$ is uniformly sampled from $[0, 0.1]$ and $\bar{\signal}_i$ is the mean of  $\signal_i$. 
The weighs in the $\ell^1$-term are taken as    $w_{\ell(\lambda)} = 2^{-\ell}$ where $\ell(\lambda)$ is the index of the downsampling-step, see \cite{obmann2020deep}.

We train all networks using the Adam \cite{kingma2014adam} optimizer with the recommended parameters for $100$ iterations and use only the best parameters of these iterations. Here, the best parameters are those which give the smallest loss on the validation set. The parameters $\eta, \beta, \gamma$ are chosen empirically and we found that $\eta = 10^{-3}$ and $\beta = \gamma = 10^{-5}$ give the best results for our approach.

\subsection{Solution of sparse  aNETT}

For  minimizing the sparse aNETT functional \eqref{eq:anett}  we use a splitting approach. For that purpose we introduce the auxiliary variable  
$\xi = \encoder(\signal)$  and rewrite  \eqref{eq:anett}   as 
the following constraint optimization problem 
\begin{equation*} 
\left\{
\begin{aligned}
&\min_{\signal, \xi} && \norm{\Ko \signal - \data}_2^2 + \alpha \norm{\xi}_{1, w} + \frac{\alpha c}{2} \norm{\signal - \auto (\signal)}_2^2 \\
&\;\text{s.t.} && \encoder(\signal) = \xi \,.
\end{aligned}
\right.
\end{equation*}
Note that we have only replaced $\encoder(\signal) $ in th $\ell^1$-term
but not in the augmented term. 
To solve the above constrained version of aNETT, we use the ADMM scheme with scaled dual variable. This results in the update scheme
\begin{align} \label{eq:admm1}
\signal_{k+1} &= \argmin_\signal \norm{\Ko \signal - \data}_2^2 + \frac{\alpha c}{2} \norm{\signal -\auto (\signal)}_2^2 \\ \nonumber
& \hspace{2cm}+\frac{\rho}{2}\norm{\encoder(\signal) - \xi_k + \eta_k}_2^2
\\ \label{eq:admm2}
\xi_{k+1} &= \argmin_\xi \alpha\norm{\xi}_{1, w} + \frac{\rho}{2} \norm{\encoder(\signal_{k+1}) - \xi + \eta_k}_2^2
\\ \label{eq:admm3}
\eta_{k+1} &= \eta_k + (\encoder(\signal_{k+1}) - \xi_{k+1}) \,,
\end{align}
where $\rho>0$ is a scaling parameter. 
The strength of the splitting type iteration 
\eqref{eq:admm1}-\eqref{eq:admm3} is that the optimization problems involved in each iterative update is  simpler and easier to solve than the original sparse aNETT minimization problem \eqref{eq:anett}, which  contains  the non-differentiable  $\ell^1$-norm as well as non-linear augmented  network term.    
In fact, the  $\xi$-update can be explicitly solved by soft-thresholding. Additionally, if we take $\auto$ being differentiable, the $\signal$-update can be  solved efficiently using gradient type iterative schemes.

\begin{table}[htb]
\centering
\scalebox{0.8}{
\begin{tabular}{l|l|l|l|l|l|l}
 & $\alpha$ & $c$ & $\rho$ & {\tt outer} & {\tt inner} & {\tt stepsize} \\ \hline \hline
\SI{0}{\percent} noise & $10^{-5}$ & $10^{2}$ & $2$ & $50$ & $10$ & $5\cdot10^{-1}$ \\
\SI{5}{\percent} noise & $5\cdot10^{-4}$ & $10^{1}$ & $2$ & $100$ & $10$ & $10^{-1}$ \\ 
adversarial & $10^{-5}$ & $10^{1}$ & $2$ & $50$ & $10$ & $5\cdot10^{-1}$
\end{tabular}}
\caption{Parameter specification for algorithm \eqref{eq:admm1}-\eqref{eq:admm3}.}
\label{tab:parameters}
\end{table}

We minimize \eqref{eq:admm1}  using gradient descent with momentum parameter $0.8$. The ADMM is initialized with $\signal_0 =  \auto (\Ro \data)$, $\xi_0 = \encoder(\signal_0)$ and $\eta_0 = 0$, where $\data$ are  the given data. Here and below $ \Ro$ denotes the filtered backprojection operator.  The parameter specifications for the minimization using  \eqref{eq:admm1}-\eqref{eq:admm3}  in various scenarios are  shown in Table~\ref{tab:parameters}. All parameters were chosen empirically to give the best results. Here, {\tt outer} refers to the total ADMM iterations, {\tt stepsize} is the stepsize and {\tt inner} is the maximal number of  iterations for the $\signal$-update step \eqref{eq:admm1}. 

\begin{figure}[htb!] 
\begin{center}
\begin{subfigure}{0.5\columnwidth}
    \includegraphics[scale=0.4,trim={3.5cm 1.1cm 3cm 1cm},clip]{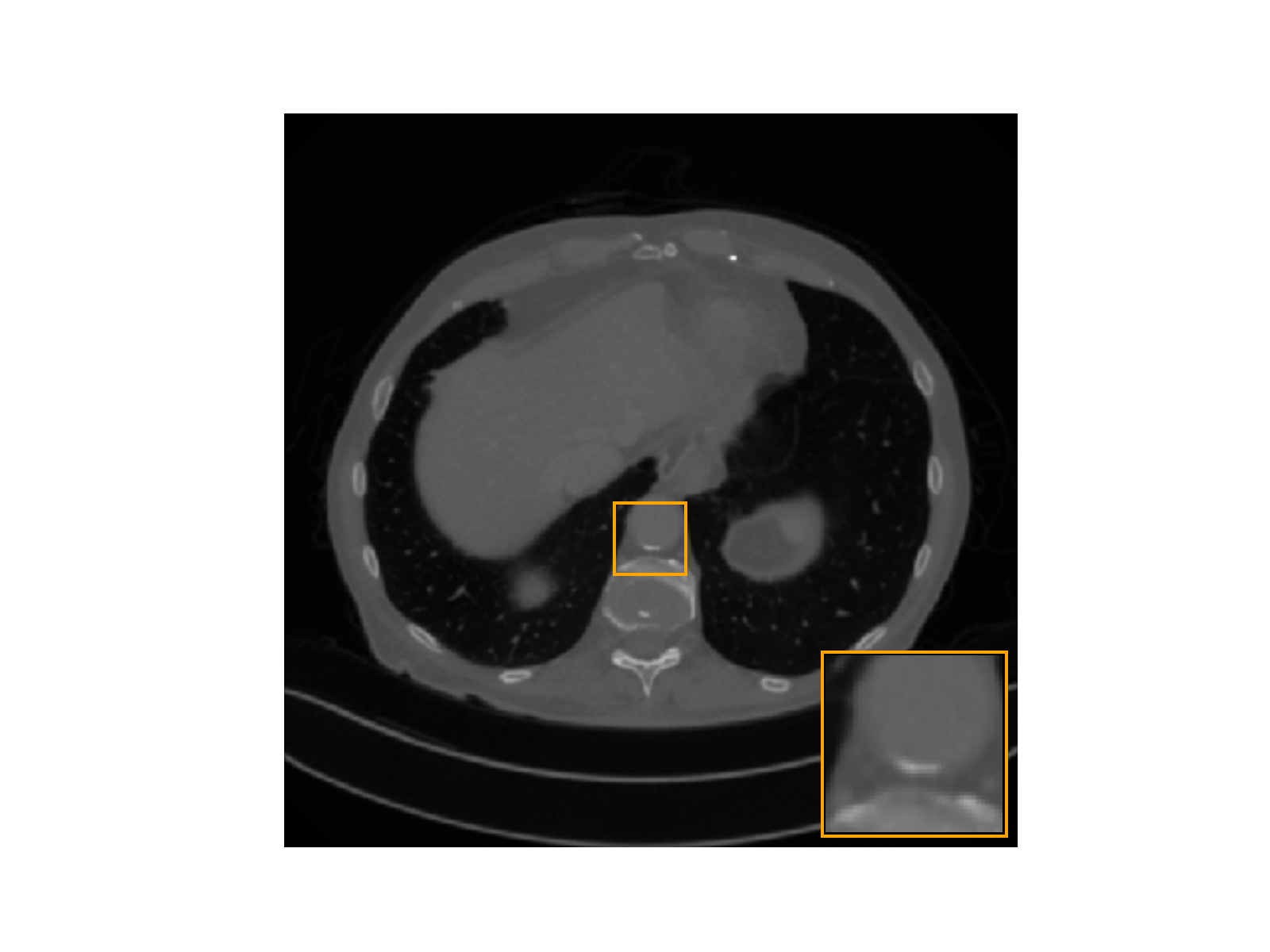}
    \vspace{-2mm}
    \caption{True}
  \end{subfigure}
  \hspace{-5mm} 
  \begin{subfigure}{0.5\columnwidth}
    \includegraphics[scale=0.4,trim={3.5cm 1.1cm 3cm 1cm},clip]{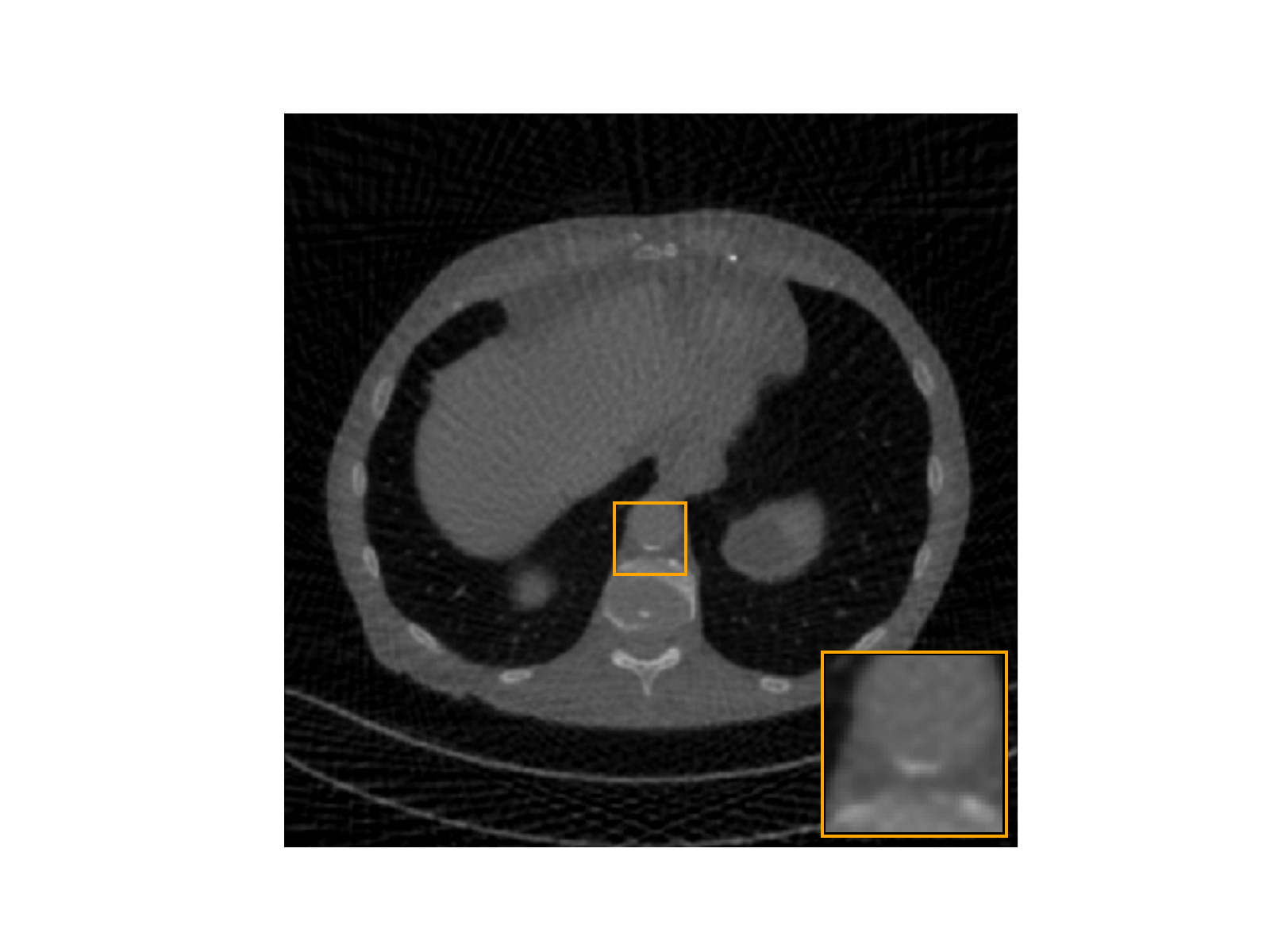}
    \vspace{-2mm}
    \caption{FBP}
  \end{subfigure} \\ \vspace{-1mm}
  \begin{subfigure}{0.5\columnwidth}
    \includegraphics[scale=0.4,trim={3.5cm 1.1cm 3cm 1cm},clip]{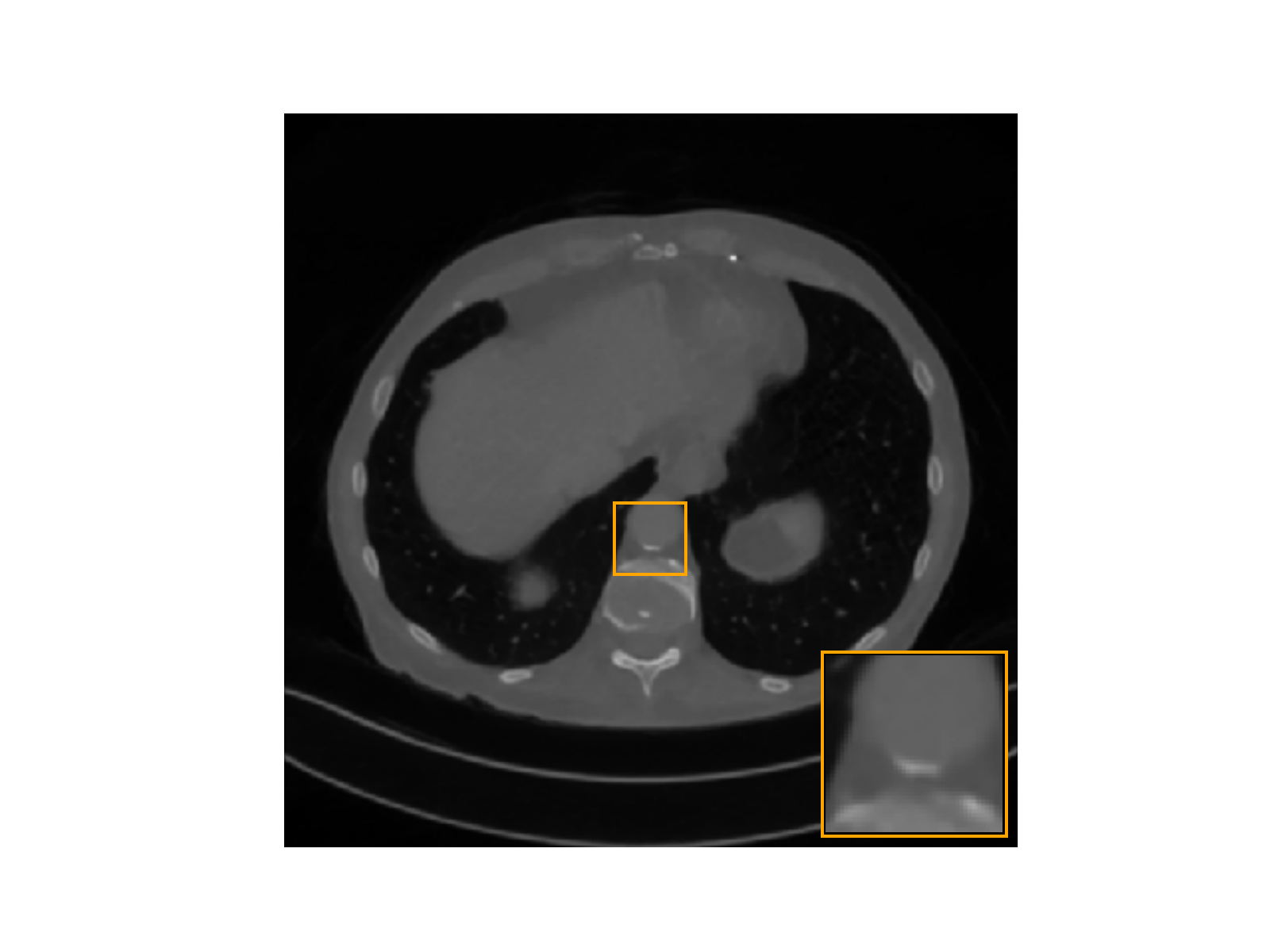}
    \vspace{-2mm}
    \caption{Network}
  \end{subfigure}
  \hspace{-5mm} 
  \begin{subfigure}{0.5\columnwidth}
    \includegraphics[scale=0.4,trim={3.5cm 1.1cm 3cm 1cm},clip]{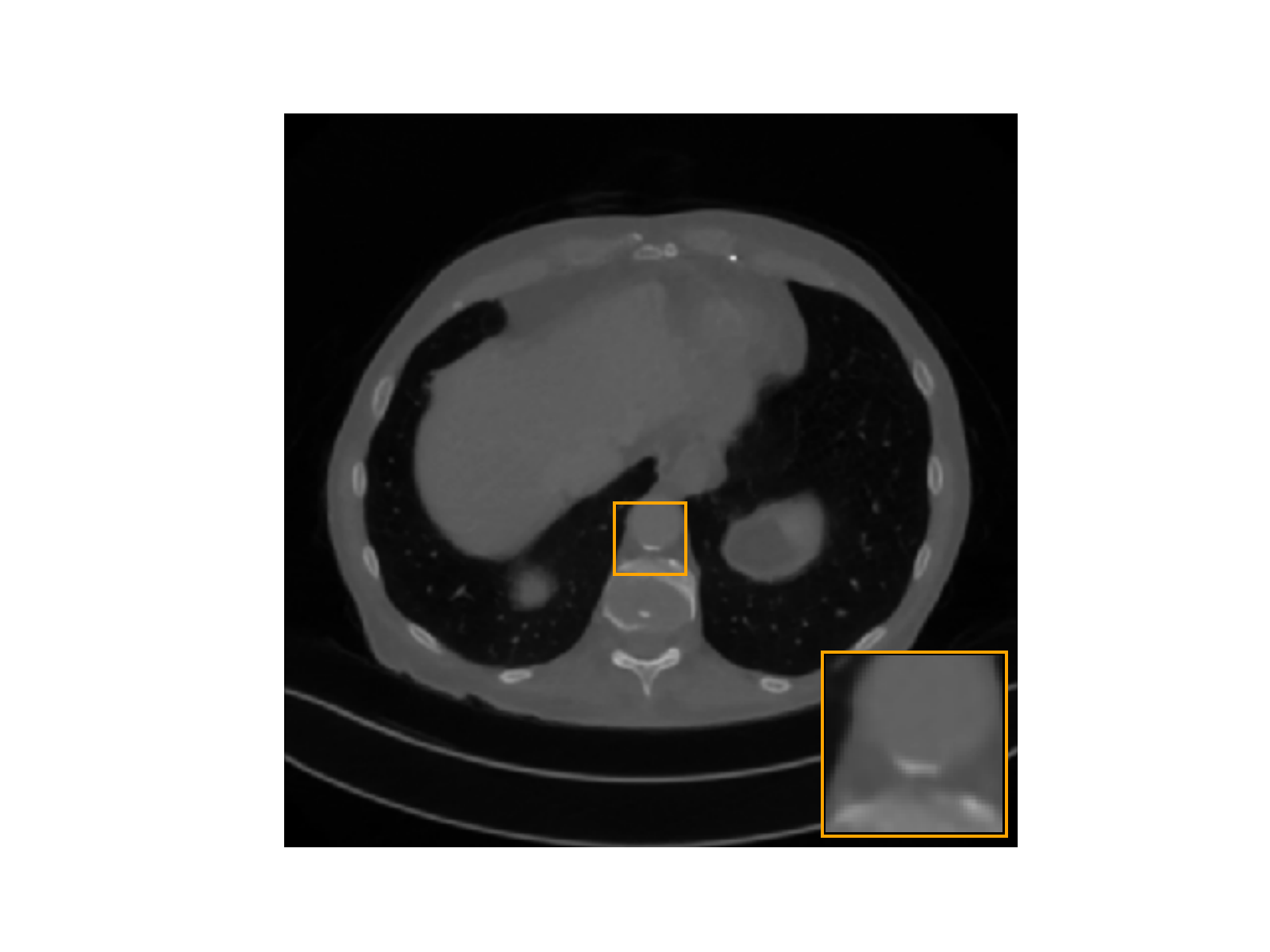}
    \vspace{-2mm}
    \caption{aNETT}
  \end{subfigure}
\end{center}
\vspace{-3mm}
  \caption{Reconstruction from simulated data.}
  \label{fig:noisefree}
\end{figure}

\subsection{Numerical results}

The first case we consider is the case of noise-free data. Figure~\ref{fig:noisefree} shows the FBP reconstruction $\signal_{\mathrm{FBP}}  = \Ro (\data) $ and the reconstruction  with the full network $\signal_{\mathrm{post}} = \auto (\signal_{\mathrm{FBP}})$ where $\auto$ is defined as above and the aNETT reconstruction  $\signal_{\mathrm{aNETT}}$. Comparing the results we see that the output of the problem adapted network $\signal_{\mathrm{post}}$ and the aNETT output $\signal_{\mathrm{aNETT}}$ are visually identical. This is because,  the test image $\signal$ is close  to the training data and therefore the considered training  procedure implies   that $\signal_{\mathrm{post}}$ is close to minimizer of the sparse aNETT.  
In comparison to the FBP we see that the aNETT was able to completely remove all the artefacts and yields an almost perfect reconstruction.

\begin{figure}[htb!] 
\begin{center}
\begin{subfigure}{0.5\columnwidth}
    \includegraphics[scale=0.4,trim={3.5cm 1.1cm 3cm 1cm},clip]{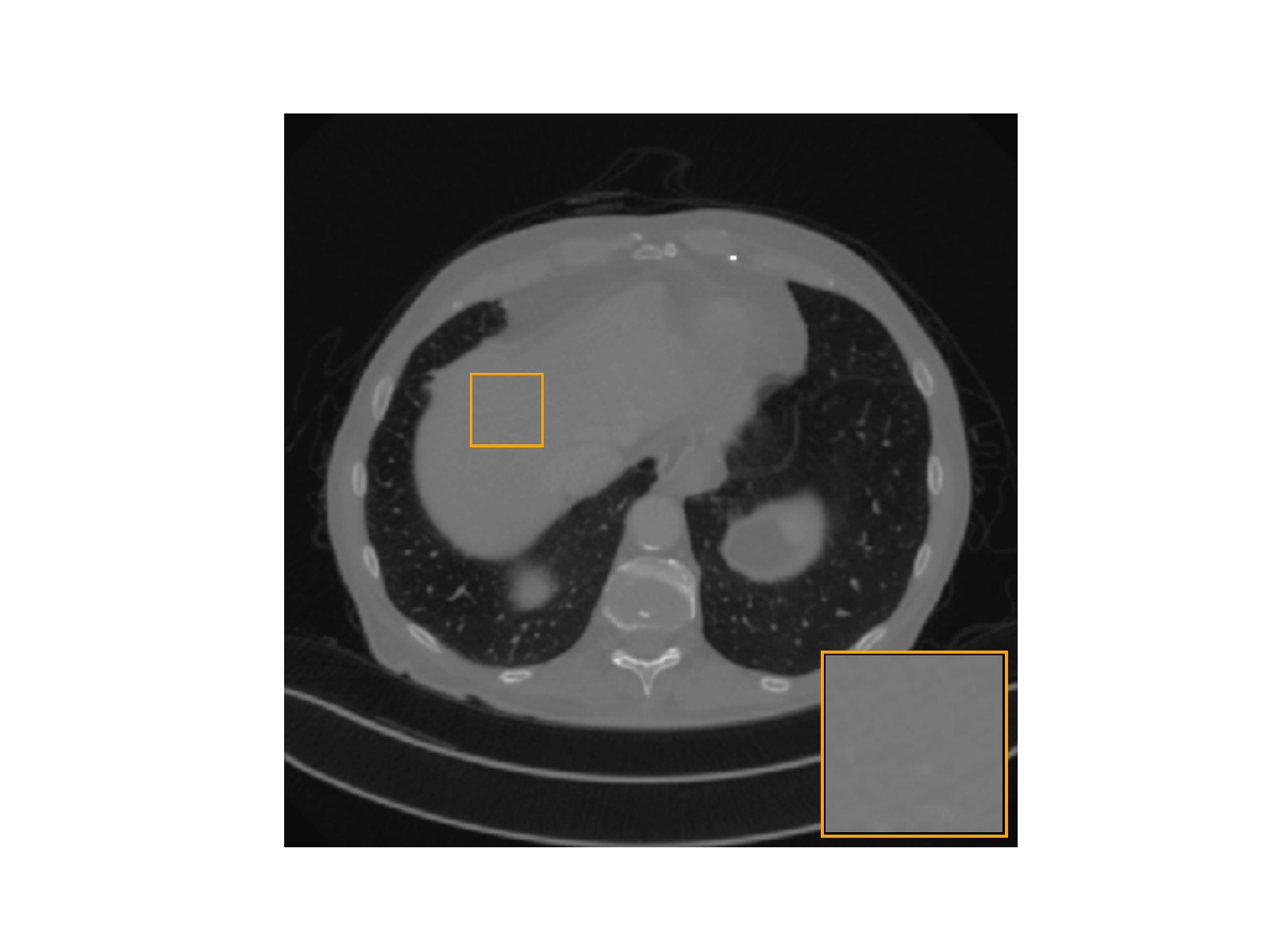}
    \vspace{-2mm}
    \caption{True}
  \end{subfigure}
  \hspace{-5mm} 
  \begin{subfigure}{0.5\columnwidth}
    \includegraphics[scale=0.4,trim={3.5cm 1.1cm 3cm 1cm},clip]{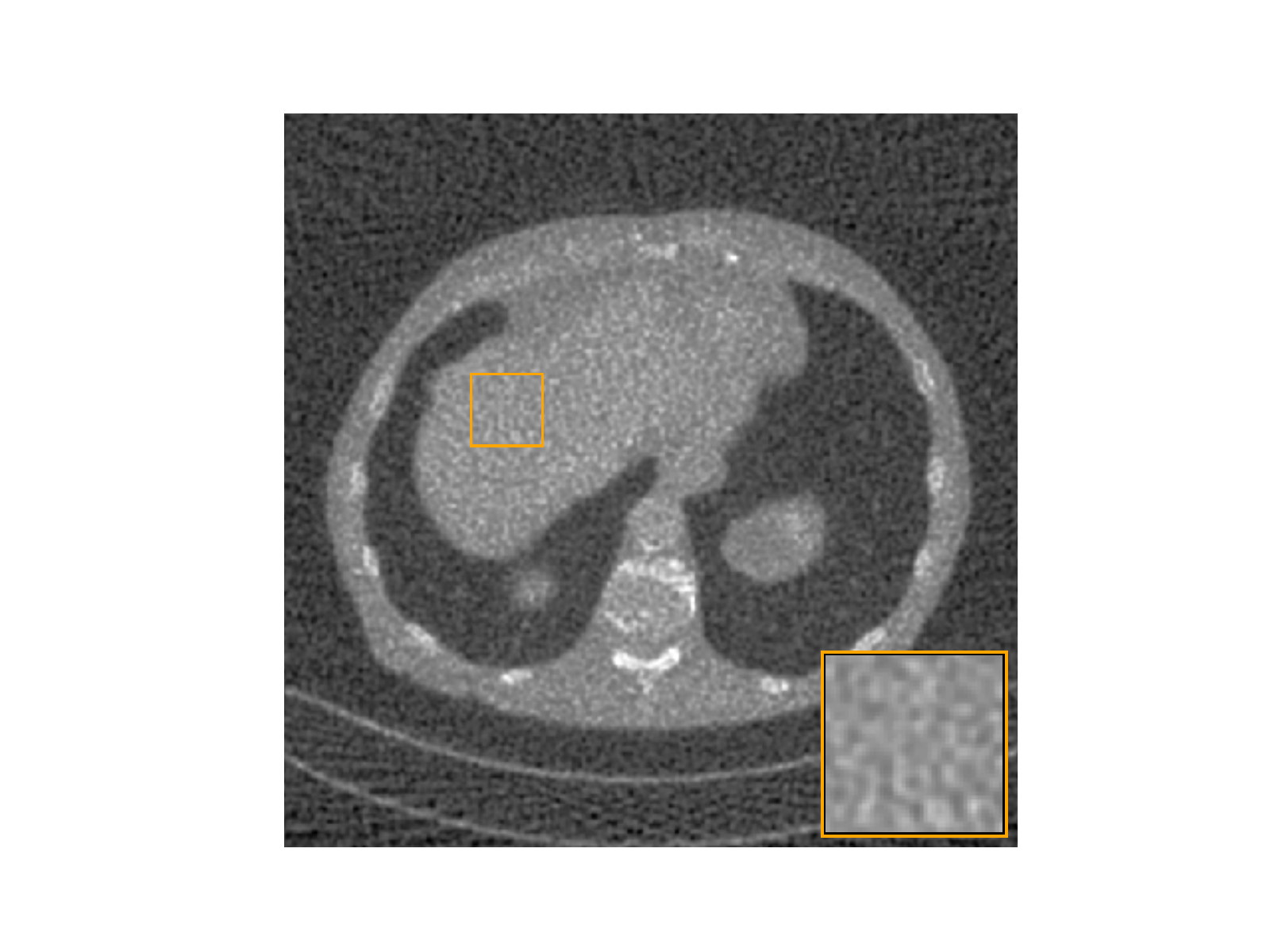}
    \vspace{-2mm}
    \caption{FBP}
  \end{subfigure} \\ \vspace{-1mm}
  \begin{subfigure}{0.5\columnwidth}
    \includegraphics[scale=0.4,trim={3.5cm 1.1cm 3cm 1cm},clip]{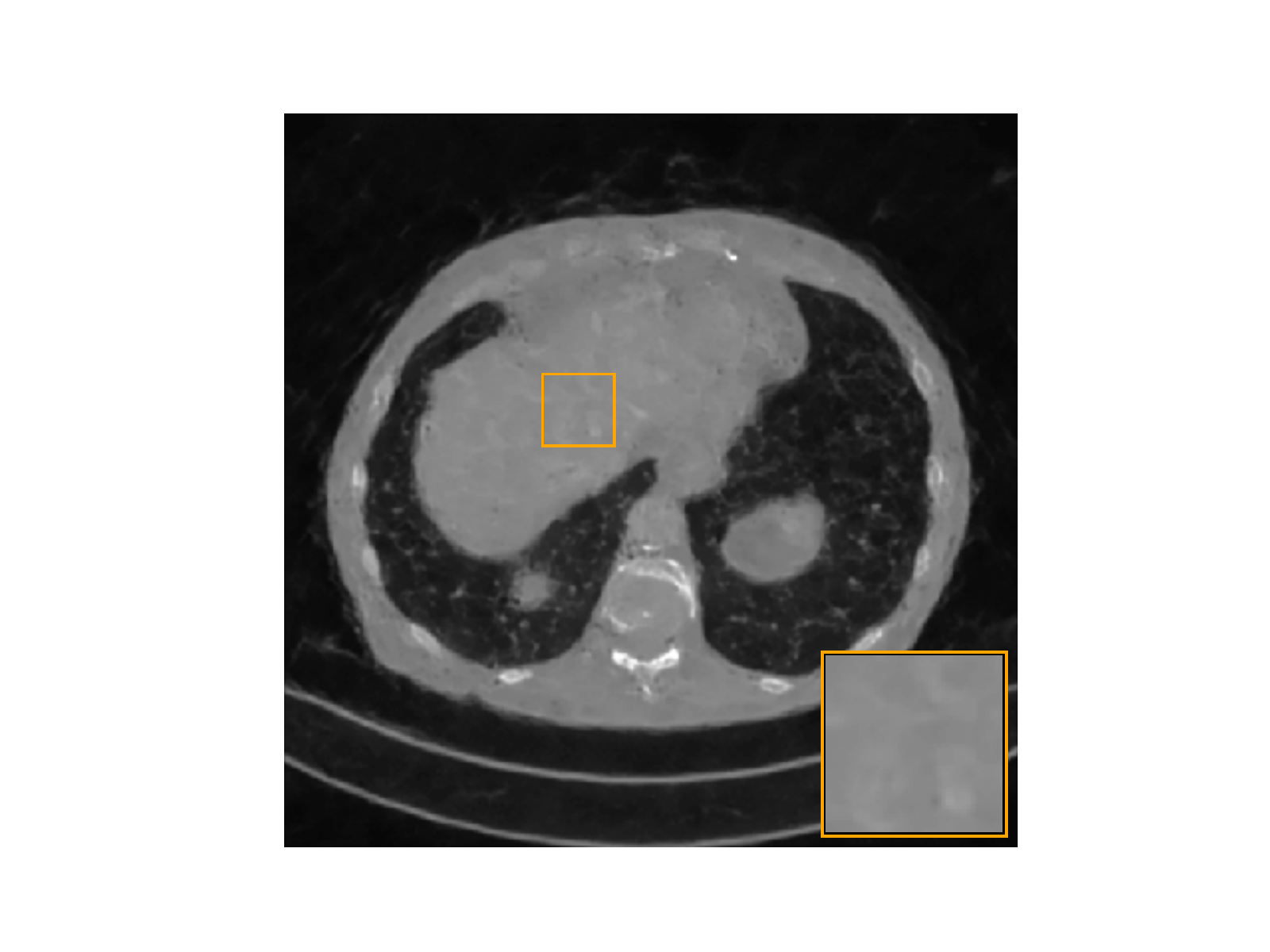}
    \vspace{-2mm}
    \caption{Network}
  \end{subfigure}
  \hspace{-5mm} 
  \begin{subfigure}{0.5\columnwidth}
    \includegraphics[scale=0.4,trim={3.5cm 1.1cm 3cm 1cm},clip]{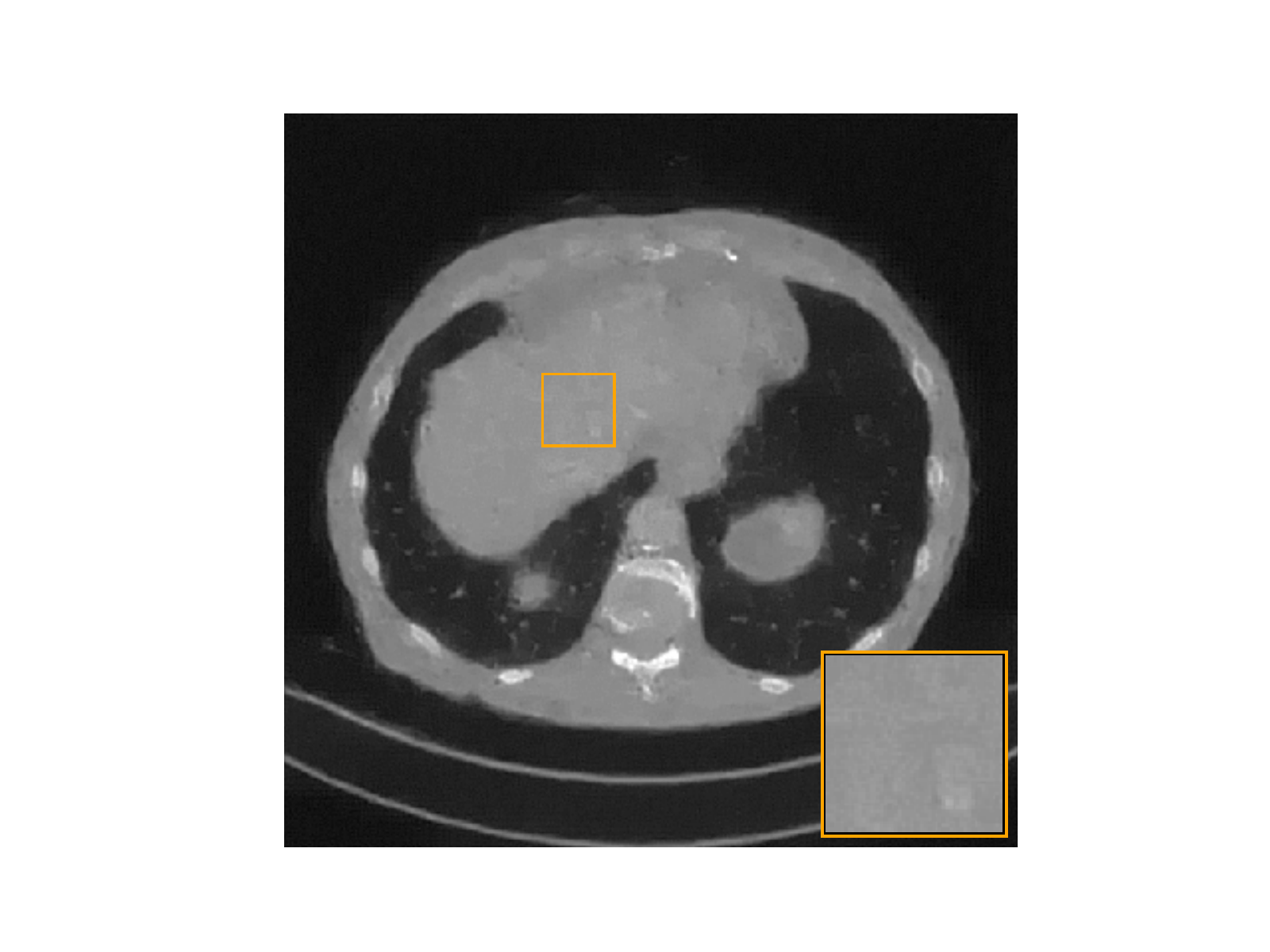}
    \vspace{-2mm}
    \caption{aNETT}
  \end{subfigure}
\end{center}
  \vspace{-3mm}
  \caption{Reconstruction from simulated data with \SI{5}{\percent} Gaussian noise. The contrast is enhanced  to emphasize the difference in the reconstructions.}
  \label{fig:noisy}
\end{figure}

To simulate noisy data we add \SI{5}{\percent} Gaussian noise to the measurement data, i.e. we use $\data_\delta = \data + 0.05 \cdot \bar{\data} \cdot \delta$ where $\bar{\data}$ is the mean of the data and $\delta$ is a standard normal distributed noise term. Reconstructions  using FBP, post-processing and the sparse aNETT  are shown in   Figure~\ref{fig:noisy}. We enhance the contrast in these images by a factor of $c = 1.7$ using the Python Pillow library \cite{clark2015pillow} to make the differences more clearly visible. 
The  post-processing  reconstruction  shows some noise-like structure on parts where the image should be mostly constant, e.g. in and around the orange square. 
We hypothesize that these noise like structures occur because the problem adapted network $\mathbf{U}$ has not been trained with noise in the data domain and hence has difficulties in reconstructing these. While we could add this to the training the networks would then likely fail on different noise models, e.g. Poisson noise.
Comparing this to the aNETT we see that this noise-like structure has been greatly reduced and we have to rely more on the sparsifying term of the regularization method to get noise-free reconstructions.

\subsection{Robustness to adversarial attack}

\begin{figure}[htb!] 
\begin{center}
\begin{subfigure}{0.5\columnwidth}
    \includegraphics[scale=0.4,trim={3.5cm 1.1cm 3cm 1cm},clip]{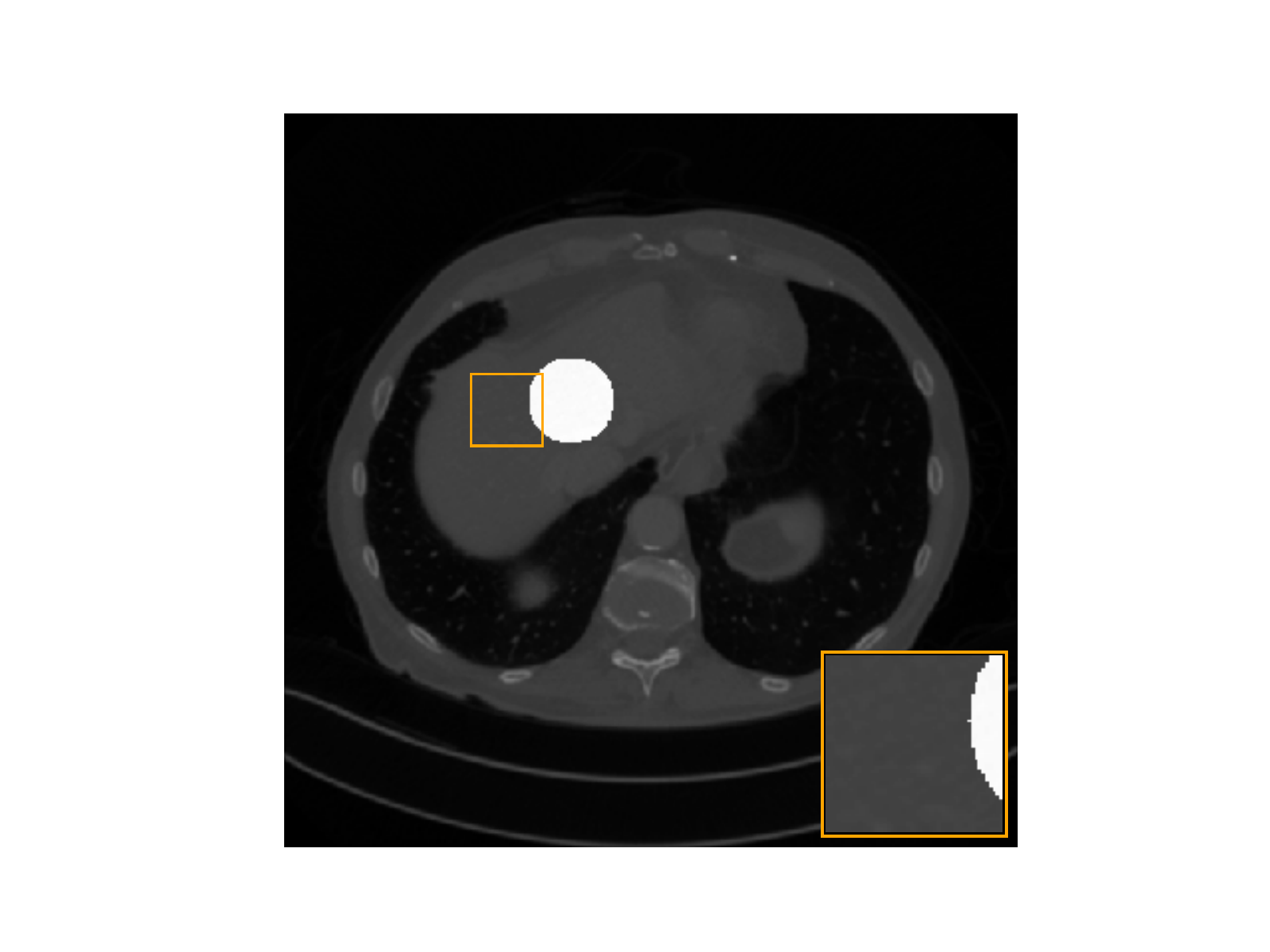}
    \vspace{-2mm}
    \caption{True}
  \end{subfigure}
  \hspace{-5mm} 
  \begin{subfigure}{0.5\columnwidth}
    \includegraphics[scale=0.4,trim={3.5cm 1.1cm 3cm 1cm},clip]{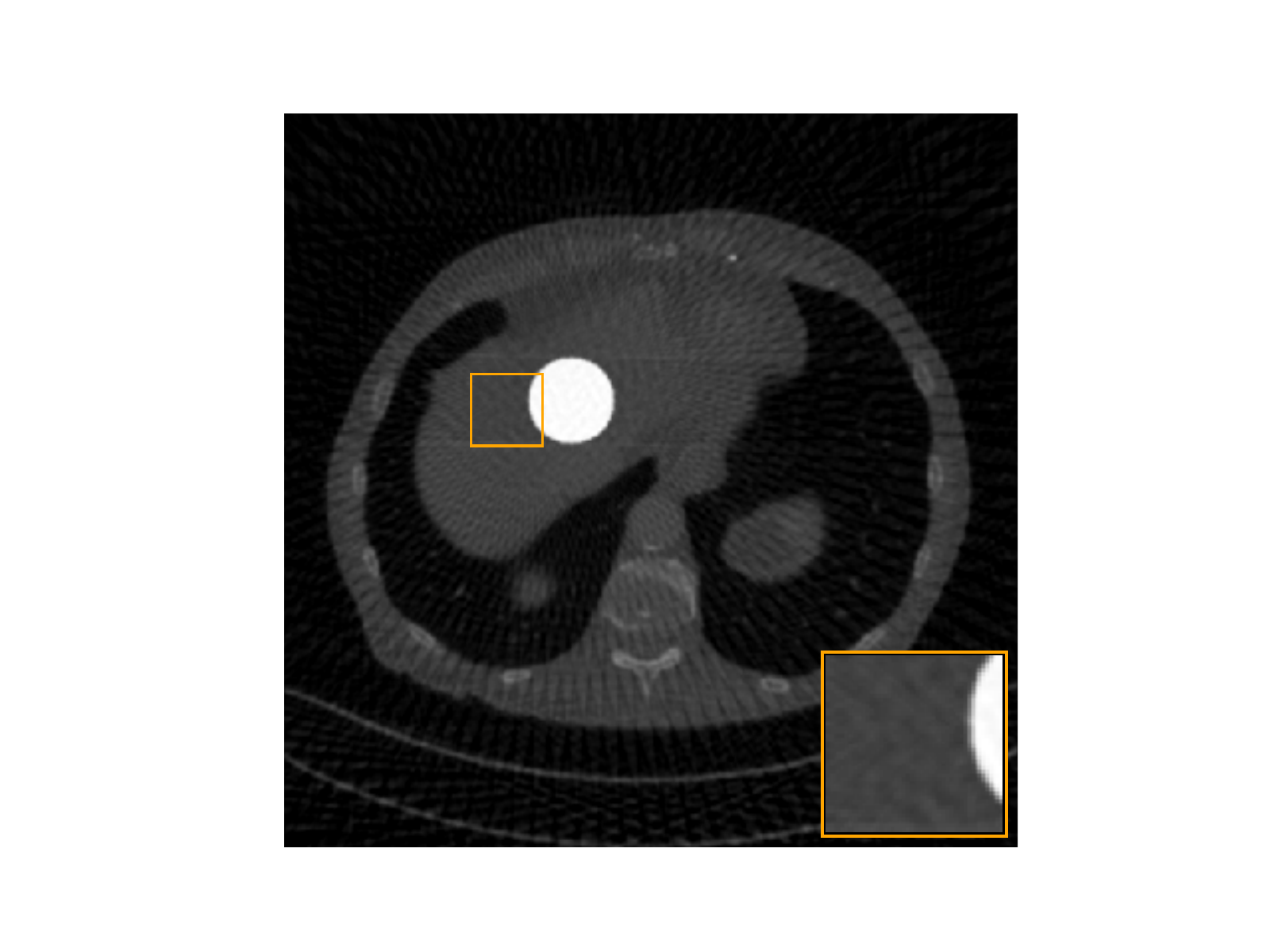}
    \vspace{-2mm}
    \caption{FBP}
  \end{subfigure} \\ \vspace{-1mm}
  \begin{subfigure}{0.5\columnwidth}
    \includegraphics[scale=0.4,trim={3.5cm 1.1cm 3cm 1cm},clip]{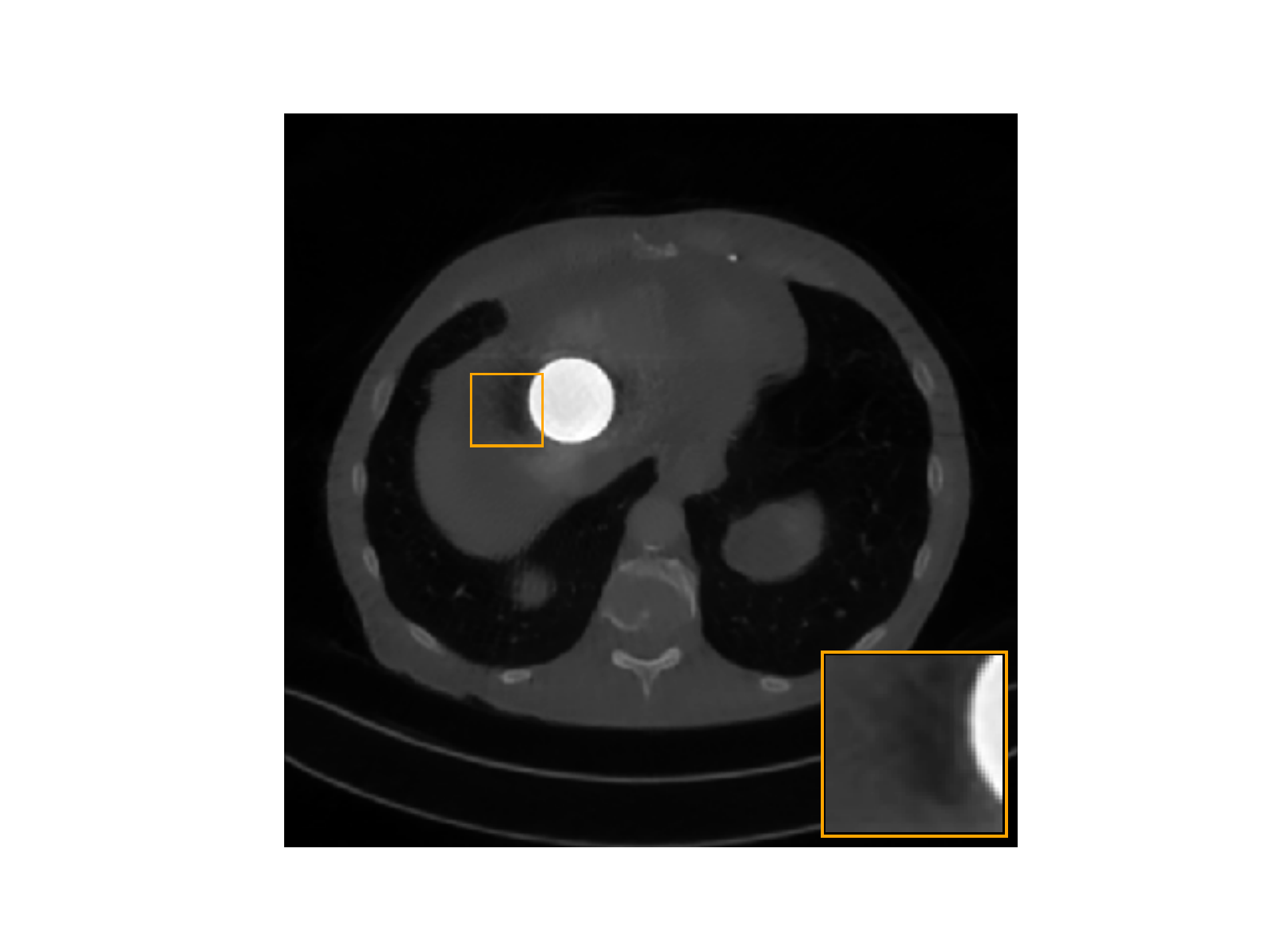}
    \vspace{-2mm}
    \caption{Network}
  \end{subfigure}
  \hspace{-5mm} 
  \begin{subfigure}{0.5\columnwidth}
    \includegraphics[scale=0.4,trim={3.5cm 1.1cm 3cm 1cm},clip]{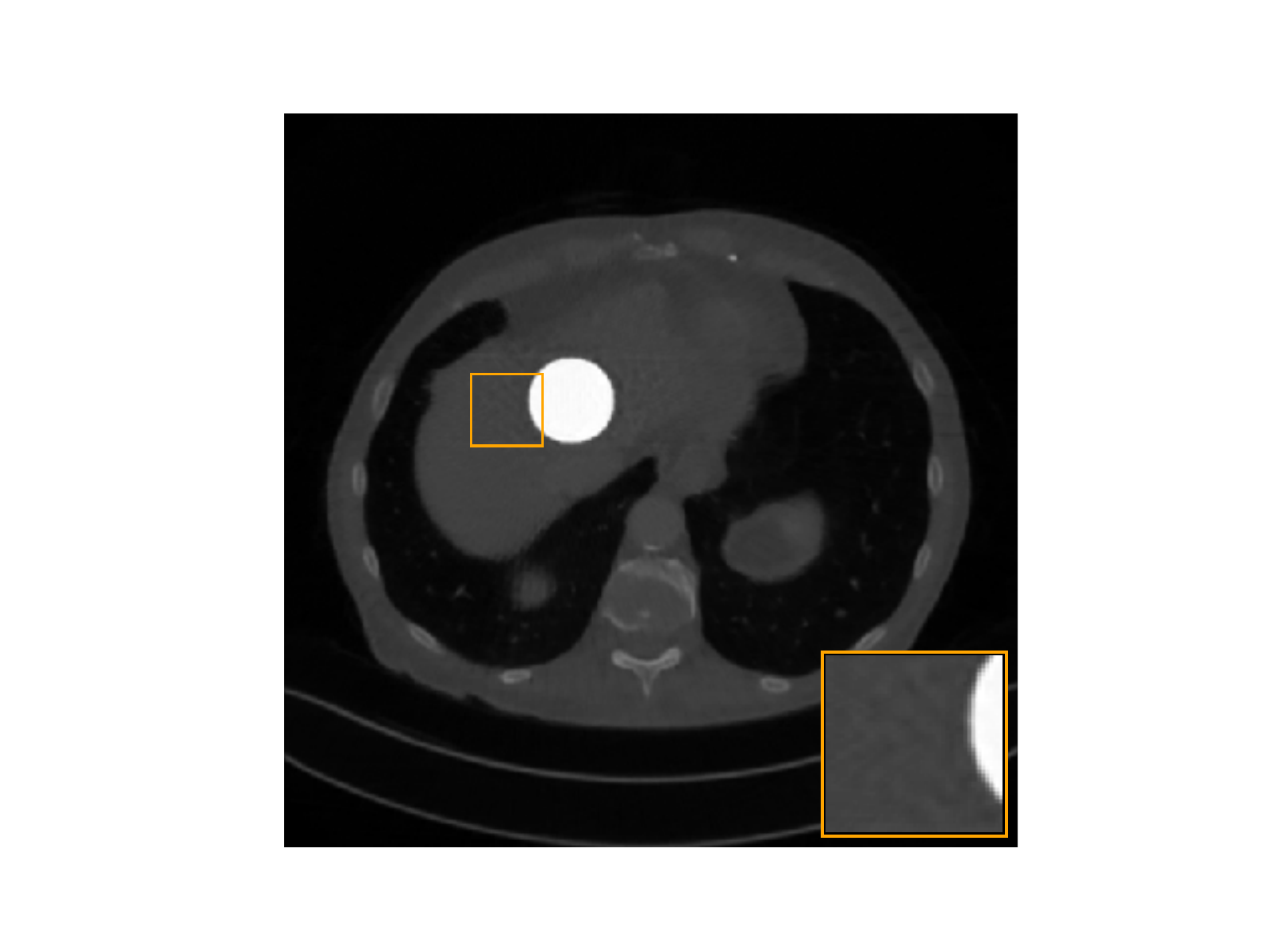}
    \vspace{-2mm}
    \caption{aNETT}
  \end{subfigure}
\end{center}
  \vspace{-3mm}
  \caption{Reconstruction from  data with additional structure.}
  \label{fig:adversarial}
\end{figure}

One particular advantage  of aNETT  over post-processing is the increased  robustness with respect to the type of image to be reconstructed.   To highlight this advantage,  as illustrated in the top left image in  Figure~\ref{fig:adversarial} we add a high intensity disc  to the  CT image shown in \ref{fig:noisefree}. The  disc represents a clear  low complexity  structure and its accurate reconstruction should be easily possible.


Figure~\ref{fig:adversarial} shows the reconstructions using the FBP, the post-processing network and the aNETT. Taking a look at the zoomed in square in these images we see that  FBP well reconstructs the circle. The post-processing  network output, however, has some dark spots close to the circle and generally shows data-inconsistent behaviour around the circle. On the other hand, using the aNETT we see that these problems do not occur. This improved accuracy  is because aNETT  takes into account the given data even for images different form the training data.

\section{Discussion}

In this paper we introduced the sparse aNETT which is a sparse reconstruction framework using a  learned regularization term and founded on a solid mathematical fundament.  As we have shown in our numerical experiments, the aNETT shows results similar to a post processing network in the case of  noise-free data phantoms  close to the training data. However, thanks to included data consistency,   the aNETT approach can much better deal with unseen phantom structure. While the chosen simple example might look  artificial, it suggests  that similar effects occur for more complex structures  in a real scenario. When considering the case of noisy data, the aNETT is able to leverage the sparsifying term   and increase robustness  with respect to noise.

While the  aNETT gives an overall more robust and stable reconstruction method, there is currently one major downside. Namely, our proposed approach relies on an iterative minimization scheme and is therefore  substantially slower than the reconstruction by a post-processing network.
Therefore the design of numerical schemes   for minimizing the sparse aNETT   functional  is a main step of future research.  Further, comparisons with different reconstruction  methods including network cascades \cite{kofler2018u,schlemper2017deep}, variational and iterative networks \cite{sun2016deep,adler2017solving,kobler2017variational} and null space networks \cite{schwab2019deep} in future work.

\bibliographystyle{IEEEbib}
\bibliography{refs}

\end{document}